\documentclass[reqno]{amsart}

\theoremstyle{definition}

\usepackage[T1]{fontenc}
\usepackage[utf8]{inputenc}
\usepackage{lmodern}
\usepackage{textcomp}
\usepackage{microtype}
\usepackage{tikz}
\usepackage{mathtools, bm, nccmath}

\usepackage[english]{babel}
\usepackage{csquotes}
\usepackage[export]{adjustbox}
\usepackage{circuitikz}

\usepackage{tabularx}

\usepackage{graphicx}
\usepackage{setspace}

\usepackage{amsmath}
\usepackage{amsfonts}
\usepackage{amssymb}
\usepackage{amsthm}
\usepackage{xcolor}%
\usepackage{soul}
\usepackage{graphicx} 
\usepackage{lipsum}
\numberwithin{equation}{section}

\newcommand*{\N}{\mathbb{N}}
\newcommand*{\Z}{\mathbb{Z}}

\newcommand*{\R}{\mathbb{R}}

\newcommand*{\V}{\mathbf{V}}
\newcommand*{\ta}{\tilde{a}}
\newcommand*{\ka}{\mathbf{k}}
\newcommand*{\taa}{\mathbf{t}}
\newcommand*{\m}{\mathbf{m}}
\newcommand*{\xa}{\mathbf{x}}

\newcommand*{\tils}{\tilde{\mathbf{S}}}
\newcommand*{\ra}{\mathbf{R}}
\newcommand*{\tilr}{\tilde{\mathbf{R}}}
\newcommand*{\tilga}{\tilde{\gamma}}

\DeclareMathOperator{\vol}{Vol}

\usepackage{thmtools}

\declaretheorem[
	name=Theorem
	]{thm}
 
\declaretheorem[
	name=Lemma,
	sibling=thm,
	]{lem}

\declaretheorem[
	name=Definition,
	style=definition,
	]{defin}

\declaretheorem[
	name=Remark,
	style=remark,
	numbered=no
	]{rem}

\newcounter{thmletter}

\declaretheorem[
    name=Theorem,
    numberlike=thmletter
]{thmA}

\newcommand{\bal}{\boldsymbol{\alpha}}

\usepackage[pdftex,pdfpagelabels]{hyperref}

\numberwithin{equation}{section}

\newcommand\restr[2]{{
  \left.\kern-\nulldelimiterspace 
  #1 
  \littletaller 
  \right|_{#2} 
  }}

\newcommand{\littletaller}{\mathchoice{\vphantom{\big|}}{}{}{}}

\makeatletter
\@namedef{subjclassname@2020}{\textup{2020} Mathematics Subject Classification}
\makeatother

\calclayout

\begin{document}
\title[Dispersion of dilated lacunary sequences]
{On the maximal volume of empty convex bodies amidst multivariate dilates of a lacunary integer sequence}

\author[E. Stefanescu]{Eduard Stefanescu}
\address{Institut f\"ur Analysis und Zahlentheorie, TU Graz, Steyrergasse 30, 8010 Graz, Austria}
\email{\href{mailto:eduard.stefanescu@tugraz.at}{eduard.stefanescu@tugraz.at}}


\begin{abstract}
Let \((a_n)_{n \in \mathbb{N}}\) be a lacunary sequence of integers satisfying the Hadamard gap condition. For any fixed dimension $d \geq 1$, we establish asymptotic upper bounds for the maximal gap in the set of dilates \(\{\bal a_n \}_{n \leq N}\) modulo 1 as  $N \to \infty$, for Lebesgue--almost all dilation vectors $\bal \in [0,1]^d$. More precisely, we prove that for any lacunary \((a_n)_{n \in \mathbb{N}}\) and Lebesgue--almost all  $\bal$, every convex set in $[0,1]^d$ of volume at least $(\log N)^{2}/N$ must contain an element of  the set \(\{\bal a_n \}_{n \leq N}\) mod 1, for all sufficiently large $N$. We also establish a generalized version of this result, where the $d$-dimensional Lebesgue measure is replaced by a general measure satisfying a certain Fourier decay condition. Our result is optimal up to logarithmic factors, and recovers as a special case a recent result for dimension $d=1$. 
\end{abstract}

\maketitle

\section{Introduction}
In this paper, we study the dispersion of multi-dimensional sequences which arise as dilates of a lacunary sequence of integers, taken modulo 1. The concept of ``dispersion'' generalizes the notion of the largest gap of a one-dimensional point set in $[0,1]$ to the multi-dimensional case of a set in $[0,1]^d$, where instead the ``largest empty range'' (under certain restrictions on the shape of admissible ranges) is studied. Of course there is no unique way for such a generalization to the multi-dimensional case, and one classical version of dispersion studies the largest volume of an empty axis-parallel box amidst a point configuration in $[0,1]^d$. In the present paper, we will use a stronger notion of dispersion, which asks for the existence of a ``large'' (in terms of volume) convex set which does not contain an element of a point set. Clearly, for any set of $N$ points in $[0,1]^d$, by the pigeon hole principle there always is an empty box of volume at least $1 / (N+1)$, and it is known that for a random point set of $N$ points in $[0,1]^d$ the probability of finding an empty box of volume larger than $c_d \log N / N$ is very small (for suitable $c_d$), see \cite{hkkr}. In this sense, allowing only empty sets of relatively small volume can be seen as a pseudo-randomness property of a point configuration, thus placing dispersion among other classical measures of pseudo-randomness of $[0,1]^d$-point sets such as discrepancy and diaphony; see \cite{dts} for a comprehensive exposition of these topics. It is clear that small dispersion can also be a very desirable property in the context of the construction of multi-dimensional point sets for sampling purposes, and accordingly questions on the existence and construction of low-dispersion point sets have received wide attention in the numerical analysis, complexity theory, and computational geometry communities; see for example the groundbreaking work of Rote and Tichy \cite{rotetichy},  and the recent works \cite{ahr,al,dj,litvak, tem, tvv}.

In this paper we will be concerned with dilated lacunary sequences \(\{\bal a_n \}_{n \in \mathbb{N}}\) modulo 1, where $(a_n)_{n \in \mathbb{N}}$ is a sequence of integers satisfying the Hadamard gap condition $a_{n+1} \geq r a_n$ for all $n \geq 1$ and some growth factor $r>1$, and $\bal$ is a dilation parameter in $[0,1]^d$ which is understood to belong to a generic set of ``typical'' parameters (with respect to some measure). Such dilated lacunary sequences have a long and rich history in analysis, probability theory, and number theory, and are in particular well-known for often exhibiting behavior which is in accordance with the typical behavior of independent random sequences. For a comprehensive treatment of this subject, we refer to the recent survey article \cite{abt}. The behavior of dilated lacunary sequences (for almost all dilation parameters) has been studied in great detail with respect to the discrepancy (see in particular Philipp \cite{phil}, Fukuyama \cite{fuku}, and Technau and Zafeiropoulos \cite{tz}) and with respect to pair correlation and gap distribution (see Rudnick and Zaharescu \cite{rudza}, Chaubey and Yesha \cite{chye}, and Yesha \cite{yesha}). Recently, Chow and Technau \cite{chow2023dispersion} investigated the size of the maximal gap for \(\{\alpha a_n \}_{n \in \mathbb{N}}\) modulo 1, and they proved that for Lebesgue-almost all $\alpha \in [0,1]$ the maximal gap is of order at most $(\log N)^{3 + \varepsilon} / N$ for all sufficiently large $N$. A crucial aspect of their work was that their method also allowed them to prove an analogous result with respect to more general measures, subject to a certain Fourier decay condition on the measure. This was particularly interesting since the range of admissible measures included the so-called Kaufman measure, a certain measure which is supported on the set of badly approximable numbers, which allowed Chow and Technau to apply their metric dispersion result in the context of Littlewood's problem, a central open problem in multiplicative Diophantine approximation. We refer to \cite{chow2023dispersion} for the details, and to \cite{bugi} for the wider context. In \cite{stefanescu2024dispersiondilatedlacunarysequences}, the author of the present paper managed to improve Chow and Technau's result by a logarithmic factor, and proved that for Lebesgue-almost all dilation factors $\alpha$, the maximal gap of \(\{\alpha a_n \}_{n \in \mathbb{N}}\) modulo 1 is of size at most $(\log N)^{2 + \varepsilon} / N$ for all sufficiently large $N$. Again, the method allowed a generalization to general measures, and an application in multiplicative Diophantine approximation. 

The purpose of this paper is to establish a multidimensional analogue of the main result from \cite{stefanescu2024dispersiondilatedlacunarysequences}, while also strengthening the bound by eliminating the $\varepsilon$-dependence in the exponent. As we will see, the multidimensional results is as strong as the one-dimensional one, even if the class of test sets is much larger than in the one-dimensional case (now it is the class of convex sets in $[0,1]^d$, rather than the class of sub-intervals of $[0,1]$). Before coming to the precise statement of our results and to the proofs, we point out one further connection to a  classical problem in geometry, which also has connections to Diophantine approximation. Danzer's problem, which has been unsolved since 1965, asks whether there exist a set of finite density in $\mathbb{R}^d$ such that the insection of this set with any convex set of volume at least 1 is non-empty. This problem, the state of research, and the connections with Diophantine approximation are explained in great detail in Adiceam's paper \cite{adiceam}. See also Solomon and Weiss \cite{sw}.

\section{Main Result} In this section we give the precise statement of the main result of this paper. The notation is explained in the next section.

\begin{defin}[Lacunary sequence]
    A sequence $(a_n)_{n\in\N}\subset \R$ is called a \textit{lacunary sequence}, if for a fixed \textit{growth factor} $r>1$:
    $$a_{n+1}\ge ra_{n}$$
    for every $n\in\N.$
    Let $\mathfrak{N}\subseteq\N$ be finite. We say a set $\{a_n\}_{n\in\mathfrak{N}}:=\{a_n:n\in\mathfrak{N}\}$ is \textit{lacunary of type $r$ if $a_m/a_n\ge r$ for any $n,m\in\mathfrak{N}$ and $m>n$}.
\end{defin}

\begin{thm}\label{main}
    Let $(a_n)_{n\in\N}$ be a lacunary sequence with growth factor $r>1$.  Let a normalized measure $\mu\in\mathcal{M}([0,1]^d)$ be given and assume that its Fourier transform decays as $\left|(\mathcal{F}\mu)(\xa)\right|\ll(1+\|\xa\|)^{-\Upsilon}$ for some $\Upsilon>0$. Then, for $\mu$-almost all $\bal\in[0,1]^d$ there exists an $N_0\in\N$ such that for all $N\ge N_0$ any convex body $\mathcal{C}\subset[0,1]^d$ with volume  Vol$(\mathcal{C})\simeq(\log N)^{2}/N$ contains at least one point of the set $\{\bal a_n\}_{n\le N}$ mod $\boldsymbol{1}$.
\end{thm}

\begin{rem}
   For the assumption $\left|(\mathcal{F}\mu)(\xa)\right|\ll(1+\|\xa\|)^{-\Upsilon}$, it does not matter which norm $\|\cdot\|$ we consider, since all norms are equivalent in finite dimensions.
\end{rem}


\begin{rem}
    As a special case for $d=1$, Theorem \ref{main} implies \cite[Theorems 2.3 and 2.4]{stefanescu2024dispersiondilatedlacunarysequences} with improved bound by a factor of $\log(N)^\varepsilon$.
\end{rem}

\section{Notation and Preparations}

Throughout this paper, the dimension $d \geq 1$ is assumed to be fixed. We write \(\mathbb{Z}^d\) for the set of all integer-lattice points, \(\mathbb{R}^d\) for the \(d\)-dimensional Euclidean space and $\mathbb{T}^d:=\R^d/(2\Z^d)$ for the \(d\)-dimensional torus. The matrix $\mathbb{I}$ represents the identity matrix. The Vinogradov symbol \(\ll\) means ''less than or equal to,'' up to constants; $\simeq$ means of the same order. Vectors are always represented by Greek or lower-case Latin letters in bold font, where the same non-bold letter with indices represents their entries, e.g., \(\boldsymbol{x} = (x_1, \ldots, x_d)^T\). We denote the zero-vector as \(\mathbf{0}\). Define the unit vectors as $\boldsymbol{\delta}_j:=(0,\dots,0,1,0,\dots,0)^T$, where $1$ is at the $j-$th position. 
The cardinality of a discrete set $\mathbf{A}$ is denoted by $|\mathbf{A}|$. We write $\lambda$ for the Lebesgue measure. For a Borel-set $\mathcal{B}\subseteq\R^d$ we define $\mathcal{M}(\mathcal{B})$ as the set of non-negative Borel measures supported on a compact subset of $\mathcal{B}$. In the following, when writing a measure $\mu$, it will always be in $\mathcal{M}([0,1]^d)$. We define $e(\,\cdotp):=\exp(2\pi i\, \cdotp)$.
\begin{rem}
Whenever we state ``for sufficiently large $N$'', it is important to note that this depends at most on the dimension $d$.
\end{rem}

 We write $\mathcal{S}(\mathbb{R}^d)$ for the Schwartz-space and $C^{\infty}_{\mathbf{K}}$ for the space of arbitrarily often continuously differentiable functions with support contained in a compact set $\mathbf{K}$. Consider two vectors $\mathbf{x},\mathbf{y}\in\R^d$, then $\mathbf{x}\cdot\mathbf{y}\in\R$ denotes the Euclidean inner product and $\mathbf{x}\otimes\mathbf{y}\in\R^{d\times d}$ denotes the Kronecker product.
 For the Fourier series expansion of a function \( f \), we denote its Fourier coefficients by $\hat{f}(\ka)=1/2^d\int_{[-1,1]^d}f(\xa)e^{-\pi i \xa\cdot\ka}d\lambda(\xa) $. 
We define the Fourier transform $\mathcal{F}:\mathcal{S}(\R^d)\rightarrow \mathcal{S}(\R^d)$, as follows: $\mathcal{F}f(\xa):=\int_{\R^d}f(\boldsymbol{\xi})e( \mathbf{\xa}\cdot\boldsymbol{\xi})d\lambda(\boldsymbol{\xi})$, see \cite{fourier} for details and basic properties. We remark $C_\mathbf{K}^{\infty}\subset\mathcal{S}$ and that every $f\in\mathcal{S}$ is bounded. 
For $f\in\mathcal{S}$ a version of the Poisson summation formula:
\begin{equation}\label{pois}
   \sum_{\mathbf{n}\in\Z^d}f(\mathbf{T}(\boldsymbol{\nu}+\mathbf{n}))=\left|\det(\mathbf{T}^{-1})\right| \sum_{\ka\in\Z^d} \mathcal{F}f\left(\mathbf{T}^{-1}\mathbf{k}\right)e(\mathbf{k}\cdot\boldsymbol{\nu}), 
\end{equation}
where $\mathbf{T}\in\R^{d\times d}$ is invertible and $\boldsymbol{\nu}\in\R^d$; see e.g.\ \cite[Proposition 1.4.2]{gröchenig2001foundations} for details. 
Next we define for $1\le p<\infty$ the space $L^p([0,1]^d):=\{f\ \textnormal{measurable}: \|f\|_{L^p(\mu)}<\infty\}$ with norm $\|f\|_{L^p(\mu)}:=\left(\int_{[0,1]^d}|f|^p d\mu\right)^{1/p}$. When considering the Lebesgue measure $\lambda$ we write $\|.\|_{L^p}$ instead of $\|.\|_{L^p(\lambda)}$. The functions  $e(\ka(\, \cdotp))$, $\ka\in\Z^d$ are orthonormal in $L^2([0,1]^d)$. 
The spaces $\ell^p(\mathfrak{N})$, $\mathfrak{N}\subseteq \Z^d$ with $1\le p<\infty$ are defined similarly to the $L^p$-spaces by exchanging $\mu$ above by the counting measure. The space $\ell^\infty(\mathfrak{N})$ is the set of bounded sequences with norm 
$\|\xa\|_{\ell^\infty}:=\sup_{i\in\mathfrak{N}}\{x_i\}$.
We recall more helpful properties: If $f$ is continuous on $\mathbb{T}^d$ and has an absolutely convergent Fourier series, then $S_N:=\sum_{|\mathbf{n}|\le N}\hat{f}(\mathbf{n})e^{\pi i \mathbf{n} \cdot}$ converges uniformly to $f$.
Moreover, we then have $2^{d/2}\|\hat{f}\|_{\ell^2}=\|f\|_{L^2(\mathbb{T}^d)}$, which is known as Parseval's identity, where $L^2(\mathbb{T}^d)$ is taken with respect to Lebesgue measure.

We further recall a version of Markov's inequality: Let $f$ be measurable, then
\begin{equation}\label{Markov'slambda}
     \mu(\{\bal\in [0,1]^d\,:\,\,|f(\bal)|\geq t\}) \leq \left(\frac{\|f\|_{L^p(\mu)}}{t} \right)^p . 
\end{equation}
Finally, we need a consequence of the first Borel Cantelli Lemma: Let $f_n$ be measurable, $(X_n)_{n\in\N}\subset\R^+$ an increasing sequence and assume that 
\begin{equation}\label{bClambda}
    \sum_{n=1}^\infty \mu(\{\bal\in [0,1]^d:  X_n - f_n(\bal)  > X_n/2\} ) < \infty.
\end{equation}
Then $f_n(\bal)>0$ for $\mu$-almost all $\bal\in [0,1]^d$ and all $n\in\N$ sufficiently large. These are all standard results, that can be found in e.g. \cite{durrett2019probability}.


A crucial ingredient in our proof will be the following result of Hadwiger, which allows us to switch from the class of all convex sets to the class of rotated rectangles. The theorem asserts that any convex set in $[0,1]^d$ necessarily contains a rotated rectangle of essentially the same volume, up to a constant factor which depends only on the dimension $d$.

\begin{thmA}[Hadwiger, \cite{378850199_0010}]\label{hard}
    Let \( \mathcal{C} \) be a convex body (compact convex set) in the \( d \)-dimensional Euclidean space with volume \(\vol(\mathcal{C}) \). Then
there exists a hyperrectangle $P$,
which covers the convex body \( \mathcal{C} \), whose volume \( \vol(P) \)  satisfies the inequality

\begin{equation}\label{a}
  \vol(P) \leq d! \vol(\mathcal{C}).
\end{equation}

 Furthermore, there exists a hyperrectangle \( Q \), which is contained in the convex body \( \mathcal{C} \), such that
\begin{equation}\label{b}
   \vol(Q) \geq \left(\frac{1}{d^d}\right) \vol(\mathcal{C}).
\end{equation}
\end{thmA}

Since Hadwiger's theorem plays a key role in the  proof of our Theorem \ref{main}, and since Hadwiger's article \cite{378850199_0010} is only available in German language, for the convenience of the reader we include a translation of Hadwiger's (short) proof of  Theorem \ref{hard} in the Appendix of this paper.\footnote{When finalizing this article we became aware that Adiceam's paper \cite{adiceam} gives a reference to a theorem of John \cite{john}, which is very similar in spirit to Hadwiger's theorem, but is formulated in terms of ellipsoids rather than rectangles. John was particularly interested in the existence and the characterization of a \emph{largest} ellipsoid contained in a convex body, while Hadwiger's paper does not attempt any such characterization. John's theorem and its proof are presented on pages 13--19 of Ball's book \cite{ball}. Hadwiger's and John's theorems seem to imply one another, but with a loss in the constants depending on $d$.}

\begin{lem}\label{lema}
    We can represent every $d$-dimensional hyperrectangle with center at the origin by rotating a hyperrectangle with axis-parallel sides and same center $\tilde{d}:=d(d-1)/2$-many times by angles $\gamma\in[0,2\pi]$ along two dimensional planes spanned by pairs of unit vectors. In other words, we can represent every rotation matrix (orthogonal matrix with determinant one) as a product of $\tilde{d}$ matrices of the form 
    \begin{equation}\label{defrotmat1}
    \mathbf{R}_{j}(\gamma):=\mathbb{I}+(\cos(\gamma)-1)\mathbf{E}_j+\sin(\gamma)\mathbf{F}_j,
\end{equation}
where $\mathbf{E}_j=(\boldsymbol{\delta}_{j-1}\otimes\boldsymbol{\delta}_{j-1}+\boldsymbol{\delta}_j\otimes\boldsymbol{\delta}_j)$ and $\mathbf{F}_j=(\boldsymbol{\delta}_{j-1}\otimes\boldsymbol{\delta}_j-\boldsymbol{\delta}_j\otimes\boldsymbol{\delta}_{j-1})$, $j\ne k$ for some $j=1,\dots,d$.
\end{lem}
\begin{rem}
    This means the matrix $\mathbf{R}_j(\gamma)$ rotates a vector by the angle $\gamma$ within the plane spanned by $\boldsymbol{\delta}_{j-1}$ and $\boldsymbol{\delta}_j$.
\end{rem}
\begin{proof} 
Without loss of generality, let us only consider cubes of side lengths $2$. Let $\mathbf{Q}$ be an arbitrary cube of this type and define $\xa^j=(x^j_1,\dots,x^j_d)^T$, $j=1,\dots,d$ as side-parallel vectors of length $1$, with  $\xa^l\cdot\xa^k=0$ for $l\ne k$ by definition.
We will construct $\tilde{d}$-matrices whose matrix product, when applied to the vectors $\xa^j$ for $j = 1, \dots, d$, equals the unit vectors $\boldsymbol{\delta}_j$ for $j = 1, \dots, d$, respectively. 

The definition \eqref{defrotmat1} yields for every vector $\mathbf{v}=(v_1,\dots,v_{j-2},0,0,v_{j+1},\dots,v_d)^T$ and every $\gamma\in\R$:
\begin{equation}\label{rotinv}
    \ra_j(\gamma)\mathbf{v}=\mathbf{v}.
\end{equation}
Writing arg$(a+ib)$ for the argument of a complex number $a+ib$, we define 
$$\gamma_k^1:= 
-\text{arg}(x_{k-1}^1+i\omega_k^1) 
$$
for $k=2,\dots,d$, where $\omega^1_k=\left(\sum_{l=k}^d\left(x_l^1\right)^2\right)^{1/2}$.
This yields $\ra_2(\gamma_2^1)\cdots \ra_d(\gamma_d^1)\xa^1=\boldsymbol{\delta^1}$. Furthermore, 
$\ra_2(\gamma_2^1)\cdots \ra_d(\gamma_d^1)\xa^j\cdot\boldsymbol{\delta}_1=0,$ 
for all $j=2,\dots,d$, hence 
\begin{equation}\label{rindexpotenz}
    \ra_2(\gamma_2^1)\cdots \ra_d(\gamma_d^1)\xa^2=(0,r_2^2x_2^2,\dots.,r_d^2x_d^2),
\end{equation}
where we define $r_j^2$, $j=2,..,d$ by the above equation \eqref{rindexpotenz}. Now define 
$$\gamma_k^2:= 
-\text{arg}(r_{k-1}^2x_{k-1}^2+i\omega^2_k) 
$$
for $k=3,\dots,d$, where $\omega^2_k:=\left(\sum_{l=k}^d\left(r_l^2x_l^2\right)^2\right)^{1/2}$. This yields $\ra_3(\gamma_3^2)\cdots \ra_d(\gamma_d^2)\ra_2(\gamma_2^1)\cdots \ra_d(\gamma_d^1)\xa^2=\boldsymbol{\delta}_2$ and by equation \eqref{rotinv} we have $\ra_3(\gamma_3^2)\cdots \ra_d(\gamma_d^2)\ra_2(\gamma_2^1)\cdots \ra_d(\gamma_d^1)\xa^1=\boldsymbol{\delta}_1$. We proceed like this and observe in every step we need to construct one matrix less, which yields that we can rotate an arbitrary cube $(d-1)+(d-2)+\dots+1=(d-1)d/2$ times to obtain an cube with axis-parallel sides. Rotation matrices have determinant one, hence we can make use of its inverses, which concludes the proof.  
\end{proof}
\begin{rem}
    Recall \(\tilde{d} = \frac{d(d-1)}{2}\).
When we consider a hyperrectangle constructed by rotating an axis-parallel hyperrectangle using the rotations 
\[
\ra_d\left(\gamma_d^{d-1}\right), \quad \ra_{d-1}\left(\gamma_{d-1}^{d-2}\right), \quad \ra_d\left(\gamma_d^{d-2}\right), \quad \dots, \quad \ra_d\left(\gamma_d^1\right),
\]
consecutively, we rename the rotation matrices to 
\[
\ra_{\tilde{d}}\left(\gamma_{d}^{d-1}\right), \quad\ra_{\tilde{d}-1}\left(\gamma_{d-1}^{d-2}\right),\quad \dots, \ra_1\left(\gamma_d^1\right),
\]
respectively. We will omit the angles whenever it is contextually clear.
\end{rem}

\section{Proof of Theorem \ref{main}}
We build upon the one-dimensional argument used by the author in \cite{stefanescu2024dispersiondilatedlacunarysequences}, which refined methods from Chow and Technau \cite{chow2023dispersion}. This approach achieves an improved bound by a factor of \((\log N)\) in comparison with \cite{chow2023dispersion}, and circumvents the delicate combinatorial arguments of \cite{chow2023dispersion} which go back to the work of Rudnick and Zaharescu \cite{rudza}. This improvement is made possible by employing exponential integrals instead of \(L^p\)-moments, where \(p = p(N)\).

\begin{proof}
In the following we will define a class of $C^\infty-$functions, each of which has its support contained in a hyperrectangle of volume $(\log N)^{2}/N$, and such that together the support of these functions essentially covers the unit cube \([0, 1]^d\). We define these functions in such a way that they indicate whether a point of the sequence $\{\bal a_n\}_{n\in\N}$ is contained in a hyperrectangle, by either being zero if not or being positive otherwise, so they are to be understood as smoothed versions of indicator functions of hyperrectangles. After preparatory steps, we will show  (heuristically speaking) that the sum of such a function over all points $\{\bal a_n\}_{n\in\N}$ is indeed always positive almost surely, meaning that the supporting hyperrectangle contains at least one element of the point set. During the proof it will be visible that the result is essentially independent of the location, form or orientation of the support of the functions, and instead depends only on its volumes and the overall cardinality of the point set. 
To move from a statement about convex sets to a statement about rotated hyperrectangles, we will construct a discrete set $\bar{\mathfrak{Q}}$ of rotated hyperrectangles whose cardinality is not too large (the cardinality will be no more than $N^{3d^2}$), such that every convex body $\mathcal{C}$ of volume $C(d)(\log N)^{2}/N$ contains a hyperrectangle in $\bar{\mathfrak{Q}}$. This will imply the statement.\\


\textbf{Step 1: Preparations:} 
Let $l$ be large enough, so that $r^l > e^2$. We define a subset $\{\tilde{a}_n\}_{n\le K}$ such that $\tilde{a}_n=a_{ln\log N}$, and $K= N/(l\log N)$, so that \begin{equation}\label{ta1}
    \tilde{a}_{n+1}\ge r^{l\log N}\tilde{a}_n\ge N^{\xi}\tilde{a_n},
\end{equation}
 for $\xi:=l\log(r)>2$, and accordingly for $m\le n$
 \begin{equation}\label{ta2}
     \tilde{a}_{m}\le N^{-\xi(n-m)}\tilde{a}_n.
 \end{equation}
To keep the notation simple, we assume here that $\log N$ is an integer and that $N$ is indeed divisible by $l\log N$. Otherwise one would need to use Gauß-brackets, which would add some small error terms without affecting the overall result.

We observe that for any fixed constant $c>0$ the set \(\{\tilde{a}_n\}_{n \le K}\) is linearly independent in \([-cN^2, cN^2] \cap \mathbb{Z}\) for sufficiently large \(N\). By utilizing equation \eqref{ta2} and a geometric series argument, this linear independence is established by the inequality:
    \begin{align}\label{estimate2}
\begin{split}
   &\sum_{j = 1}^{K_0 - 1} cN^2 \tilde{a_j} \le cN^2\sum_{j=1}^{K_0-1}N^{-\xi(K_0-j)}\ta_{K_0} \\
   =& cN^{2-\xi K_0}\frac{N^{\xi(K_0-1)}-1}{N^{\xi}-1}N^{\xi}\ta_{K_0}\le 2cN^{2-\xi}\ta_{K_0}<\tilde{a}_{K_0},
    \end{split}
\end{align} 
for every \(K_0 \in [2, K] \cap \mathbb{Z}\) and $N$ sufficiently large, and since $\xi>2$.
In the sequel we assume without loss of generality, \(l = 1\), since all arguments can be directly adapted to the case \(l \ne 1\).

We define the following quantities
\begin{equation}\label{PQR}
    Q:=1000d^22^d\log N \max\{\Upsilon^{-1},1\},\ M:=Q^2,
\end{equation}
and the set of scaling matrices, translation vectors and rotation matrices 
\begin{align}\label{matricesSRT}
    \begin{split}
        \mathfrak{S}&:=\left\{\mathbf{S}\in\R^{d\times d}\middle|\mathbf{S}=\begin{pmatrix}
s_1 & 0 & \cdots & 0 \\
0 & s_2  & \cdots & 0 \\
\vdots & \vdots  & \ddots & \vdots \\
0 & 0 &  \cdots & s_d
\end{pmatrix}, \prod_{i=1}^{d}s_i=2^{-d}\frac{M}{N} \right\},\\
\boldsymbol{\tau}&:=\left\{\mathbf{t}\in\R^d\middle|\mathbf{t}\in[0,1]^d\right\},\\
\mathfrak{R}&:=\left\{\mathbf{R}\in\R^{d\times d}\middle|\mathbf{R}\text{ orthogonal },\det(\mathbf{R})=1\right\}.
    \end{split}
\end{align}
We define

\begin{equation}
    \mathfrak{Q}:=\Bigl\{\mathbf{Q}=\mathbf{R}\mathbf{S}\cdot[-1,1]^d+\mathbf{t}\Bigm|\mathbf{R}\in\mathfrak{R},\mathbf{S}\in\mathfrak{S},\mathbf{t}\in\boldsymbol{\tau}\Bigr\}.
\end{equation}
Thus $\mathfrak{Q}$ is a set of rotated and translated hyperrectangles of fixed volume $M/N$.
     Next, define 
     $$\mathfrak{M}:=\left\{\mathbf{V}\in\R^{d\times d}\middle|\V^{-1}=\mathbf{R}\mathbf{S},\ \mathbf{R}\in\mathfrak{R},\ \mathbf{S}\in\mathfrak{S}\right\}.$$


Let \( f \in C_{[-1,1]^d}^{\infty} \) be a real, non-negative and even function such that $f\le 1$, and \(\int_{\mathbb{R}^d} f(\xa) \, d\lambda(\xa) = 1\). This implies that \(\mathcal{F}f\) (the Fourier transform of \(f\)) is real and even, \(\mathcal{F}f(0) = 1\), \(|\mathcal{F}f| \le 1\) and $\|f\|_{L^2}\le\|f\|_{L^1}=1$. Such a function exists, see e.g. \cite[Lemma 1.26]{demeter2020fourier}.
Furthermore let $\V\in\mathfrak{M}$ and $\taa\in\boldsymbol{\tau}$, hence $\mathbf{Q}=\V^{-1}\cdot[-1,1]^d+\taa\in\mathfrak{Q}$, and
define the following function:
\begin{equation}\label{functionmainres}
    \omega_{N,\mathbf{Q}}(\bal):=\sum_{\mathbf{u}\in\Z^d}\sum_{n=1}^{N/Q}f\left(\mathbf{V}(\tilde{a}_n\bal -\mathbf{z}_{\mathbf{t}}-\mathbf{u})\right).
\end{equation}
The sum over \(\mathbf{u} \in \mathbb{Z}^d\) accounts for the fractional parts in each component of the evaluated vector; the matrix \(\mathbf{V}\in\mathfrak{M}\) dilates the support of $\omega_{N,\mathbf{Q}}$, and the vector \(\mathbf{z}_{\taa}\in\R^d\) is defined such that it translates the support of $\omega_{N,\mathbf{Q}}$ to be centered in $\mathbf{t}$. This means we define $\mathbf{z}_\taa$ such that $\V \mathbf{z}_\taa=\taa$, and hence $f(\V(\tilde{a}_n\bal-\mathbf{z}_\taa-\mathbf{u})=f(\V(\ta_n\bal-\mathbf{u})-\taa)$.

Since the support of $\omega_{N,\mathbf{Q}}$ is contained in a $d$-dimensional hypercube, we can imagine $\mathbf{S}$ as an operator dilating the edges of the hypercube resulting in a hyperrectangle. By Lemma \ref{lema}, we can represent the matrix $\mathbf{R}$ as a succession of rotations of the resulting hyperrectangle by the angles $\gamma_1,\dots,\gamma_{\tilde{d}}$ sequentially within a plane spanned by two unit vectors, where $\tilde{d}:=d(d-1)/2$. 

Using Poisson summation formula \eqref{pois}, and noting that $\det(\V^{-1})=2^{-d}M/N$ we get:
$$\omega_{N,\mathbf{Q}}(\bal)=2^{-d}\frac{M}{N}\sum_{\ka\in\Z^d}\sum_{n=1}^{N/Q}(\mathcal{F}f)\left(\V^{-1}\mathbf{k}\right)e(\ka\cdot(\tilde{a}_n \bal -\mathbf{z}_{\mathbf{t}})).$$
We now define:
\begin{align}\label{minusQ}
    \begin{split}
       \tilde{\omega}_{N,\mathbf{Q}}(\bal)&:=2^{-d}\frac{M}{N}\sum_{\ka\in\Z^d}\sum_{n=1}^{N/Q}(\mathcal{F}f)\left(\V^{-1}\mathbf{k}\right)e(\ka\cdot(\tilde{a}_n \bal -\mathbf{z}_{\mathbf{t}}))-2^{-d}Q\\\
&=2^{-d}\frac{M}{N}\sum_{\ka\in\Z^d\backslash\{\mathbf{0}\}}\sum_{n=1}^{N/Q}(\mathcal{F}f)\left(\V^{-1}\mathbf{k}\right)e(\ka\cdot(\tilde{a}_n \bal -\mathbf{z}_{\mathbf{t}})).
    \end{split}
\end{align}
Here we used that $(\mathcal{F}f)(0)=1$ and $M/Q=Q$. 

Since after rotating the support of $f$ by $\mathbf{R}$ the total number of hyperrectangles needed to cover $[0,1]^d$ via translations does not change, up to a finite factor dependent on the dimension $d$, we can assume without loss of generality $\mathbf{R}=\mathbb{I}$, thus $$\V^{-1}=\mathbf{S}=
\begin{pmatrix}
s_1 & 0 & \cdots & 0 \\
0 & s_2 & \cdots & 0 \\
\vdots & \vdots & \ddots & \vdots \\
0 & 0 & \cdots & s_d
\end{pmatrix}.
$$
We truncate the first sum in equation \eqref{minusQ}, by summing $\ka$ over 
$$\mathbf{\Sigma}:= \left(\bigtimes_{i=1}^d[-N^{\frac{1}{d}}/s_i-1,N^{\frac{1}{d}}/s_i+1]\cap\Z^d\right)\backslash\{\mathbf{0}\},$$
where $|\mathbf{\Sigma}|\le \prod_{i=1}^d(2 N^{\frac{1}{d}}/s_i+3)\le(2^{2d+1})N^2/M$ for $N$ large enough. We work with
       $$\omega^{*}_{N,\mathbf{Q}}(\bal):=2^{-d}\frac{M}{N}\sum_{\ka\in\mathbf{\Sigma}}\sum_{n=1}^{N/Q}(\mathcal{F}f)\left(\V^{-1}\mathbf{k}\right)e(\ka\cdot(\tilde{a}_n \bal -\mathbf{z}_{\mathbf{t}})).$$
Let us define $\mathbf{\Sigma}':= \left(\bigtimes_{i=1}^d[-N^{\frac{1}{d}}/s_i,N^{\frac{1}{d}}/s_i]\cap\Z^d\right)\backslash\{\mathbf{0}\}$. Using the rapid decay of $\mathcal{F}f$ (which follows from $\mathcal{F}$ being an automorphism and from the definition of Schwartz functions), we obtain, for any $\Xi > 0$,
\begin{align}\label{trunc}
\begin{split}
     &\left|2^{-d}\frac{M}{N}\sum_{\ka\notin\mathbf{\Sigma}}\sum_{n=1}^{N/Q}(\mathcal{F}f)\left(\V^{-1}\mathbf{k}\right)e(\ka\cdot(\tilde{a}_n \bal -\mathbf{z}_{\mathbf{t}}))\right|\ll \left|Q\sum_{\ka\notin\mathbf{\Sigma}}\frac{1}{\left(1+|\V^{-1}k|\right)^\Xi}\right|\\
     \ll& Q\int_{\R^d\backslash\mathbf{\Sigma}'}  \frac{1}{\left(1+|\V^{-1}\xa|\right)^\Xi}d\lambda(\xa)=\frac{N}{Q}\int_{\R^d\backslash\left[-N^{\frac{1}{\delta}},N^{\frac{1}{\delta}}\right]^d}  \frac{1}{\left(1+|\mathbf{y}|\right)^\Xi}d\lambda(\mathbf{y})\le O\left(N^{-1000d}\right),
\end{split} 
\end{align}
for $N$ and $\Xi$ sufficiently large.

Consequently, we can write:

$$\tilde{\omega}_{N,\mathbf{Q}}(\bal)=\omega^{*}_{N,\mathbf{Q}}(\bal)+O\left(N^{-1000d}\right).$$
\begin{rem}
   In the paper \cite{stefanescu2024dispersiondilatedlacunarysequences} by the author, there was a minor inconsistency regarding the above estimate, which is correctly stated in the present paper. Specifically, we need to sum over essentially \( N^2 \) rather than \( N\log(N)^\varepsilon \). Theorem \ref{main} is stronger and provides correctness of Theorems 2.4-2.6 in \cite{stefanescu2024dispersiondilatedlacunarysequences}.
\end{rem}
Define 
\begin{equation}\label{pwiezuvor}
    p_{\ka}:=(\mathcal{F}f)\left(\V^{-1}\ka\right)e(-\ka \cdot\mathbf{z}_{\mathbf{t}}),
\end{equation}
where in particular $|p_{\ka}|\le1$. \\

\textbf{Step 2: Moment estimates (Lebesgue measure):}
To make the proof more accessible, we will first prove the result in the case where $\mu=\lambda$ is the Lebesgue measure. Subsequently, we will extend the argument to the more general case in the next step.
We estimate the following integral, which will be useful later, when using the first Borel--Cantelli Lemma.

\begin{align}\label{lebesguemeasure}
\begin{split}
   &\left\|e^{-\frac{1}{10}\omega^*_{N,\mathbf{Q}}}\right\|_{L^1}=\int_{[0,1]^d}\exp\left(-\frac{1}{10}\omega^*_{N,\mathbf{Q}}(\bal)\right)d\lambda(\bal)\\
=&\int_{[0,1]^d}\prod_{n=1}^{\frac{N}{Q}}\exp
\left(-\frac{1}{10\cdot2^{d}}\frac{M}{N}\sum_{\ka\in\mathbf{\Sigma}}p_{\ka}e(\ka\cdot\tilde{a}_n\bal)\right)d\lambda(\bal)\\
\le&\int_{[0,1]^d}\prod_{n=1}^{\frac{N}{Q}}\Biggl(1-\frac{1}{10\cdot2^{d}}\frac{M}{N}\sum_{\ka\in\mathbf{\Sigma}}\Bigl(p_{\ka}e(\ka\cdot\tilde{a}_n\bal)\Bigr)\\
&+\biggl(\frac{1}{10\cdot2^{d}}\frac{M}{N}\sum_{\ka\in\mathbf{\Sigma}}\bigl(p_{\ka}e( \ka\cdot\tilde{a}_n\bal)\bigr)\biggr)^2\Biggr)d\lambda(\bal).
\end{split}
\end{align}
Here, we can remove the absolute values in the first line, since \(\omega^*_{N,\mathbf{Q}}\) is real, making \(\exp\left(\frac{1}{10} \omega^*_{N,\mathbf{Q}}\right)\) positive. The transition from the second to the third line involves applying the Taylor series expansion for exponential functions. Specifically, for $y\in\R$, and \(|y| \le 1\), we have \(e^y \le 1 + y + y^2\). 


We define the function $h:=f(\V(\cdot/2-\mathbf{z}_\taa))$. Then three substitutions and the fact that \( f \) is supported in \( [-1,1]^d \) yield
$$\hat{h}(\ka)=\frac{M}{2^dN}\mathcal{F}f\left(\V^{-1}\ka\right)e(-\ka\cdot\mathbf{z}_\taa)=\frac{M}{2^dN}p_\ka.$$
Using 
that \( \mathcal{F}f \) is strictly radially decreasing, we split the sum over \( \mathbb{Z}^d \) into disjoint subsets corresponding to lattice points inside disjoint translates of the hyperrectangle \( \V\cdot[-1,1)^d \). 

Let $\Xi>0$, and $\mathbf{m}\in2\Z^d$ be arbitrary. Then 
for $\mathbf{k} \in V([-1,1]^d + \mathbf{m})$, we have 
$\V^{-1} \mathbf{k} \in [-1,1]^d + \mathbf{m}$ and we see
\[
|\mathcal{F} f (\V^{-1} \mathbf{k})| = O\big(1 / (1 + |\mathbf{m}|)^{\Xi}\big),
\]
which yields

\begin{align}
\begin{split}
&\sum_{\ka\in\Z^d}|\mathcal{F}f(\V^{-1}\ka)|=\sum_{\mathbf{m}\in2\Z^d}\sum_{\ka\in\V([-1,1)^d+\mathbf{m})\cap\Z^d}|\mathcal{F}f(\V^{-1}\ka)|\\ \le&\sum_{\mathbf{m}\in2\Z^d}O\big(1 / (1 + |\mathbf{m}|)^{\Xi}\big)\sum_{\ka\in\V([-1,1)^d+\mathbf{m})\cap\Z^d}1\simeq\frac{N}{M}\sum_{\mathbf{m}\in2\Z^d}O\big(1 / (1 + |\mathbf{m}|)^{\Xi}\big)
\end{split}
\end{align}
This, together with an almost identical argument as \eqref{trunc} yields
\[
\sum_{\ka \in \mathbf{\Sigma}} |\hat{h}(\ka)| \ll \frac{M}{2^dN}\frac{N}{M}\sum_{\mathbf{m}\in2\Z^d}O\big(1 / (1 + |\mathbf{m}|)^{\Xi}\big) + O\left(N^{-1000d}\right)<\infty.
\]
Hence $\hat{h}$ is absolutely summable independent of $N$ and $\taa$, due to the rapid decay of $\mathcal{F}f$. This implies, that
\begin{equation}
    \sum_{\ka\in\Z^d}\frac{M}{2^dN}p_\ka e(\ka\cdot\xa )=h(\xa)=f(\V(\frac{\xa}{2}-\mathbf{z}_\taa))\le1,
\end{equation}
for all $\xa\in\mathbb{T}^d$.
Utilizing the above together with a truncation argument as in \eqref{trunc} we get
\begin{equation}\label{new1}
   \left|\frac{1}{10\cdot 2^{d}} \frac{M}{N} \sum_{\ka \in \mathbf{\Sigma}} p_{\ka} e( \ka \cdot \bal \tilde{a}_n)\right| \le \left|\frac{1}{10} h(2\bal \tilde{a}_n) + O\left(N^{-1000d}\right)\right| < 1,
\end{equation}
Furthermore, we observe that \( \omega^{*}_{N,\mathbf{Q}} \) is real. This explains the estimate in \eqref{lebesguemeasure}.

Next, we rewrite:
\[
e(\ka\cdot\tilde{a}_n\bal) = \prod_{i=1}^d e(\tilde{a}_n k_i\alpha_i),
\]
By Fubini's theorem it is sufficient to regard every factor individually. Hence, we can treat this problem similarly to the one-dimensional case. 
Due to the linear independence of \(\{\tilde{a}_n\}_{n \leq K}\) within the interval \([-N, N]\), see estimate \eqref{estimate2} and the \(L^2[0,1]\)-orthonormality property, the integrals involving products of phases \(e(\tilde{a}_n k_j \alpha_j),\ j=1,\dots,d\) with other phases of significantly higher frequencies vanish, which yields the second inequality below. The inequality in the fourth row is a direct consequence of Bernoulli's and the triangle inequality. In the penultimate step, we used Parseval's identity. The final step is doing a substitution. Continuing from equation \eqref{lebesguemeasure}, we have

\begin{align}\label{lebesguemeasure2}
\begin{split}
\left\|e^{-\frac{1}{10}\omega^*_{N,\mathbf{Q}}}\right\|_{L^1}\ll&\prod_{n=1}^{\frac{N}{Q}}\Biggl(\int_{[0,1]^d}1-\frac{1}{10\cdot2^{d}}\frac{M}{N}\sum_{\ka\in\mathbf{\Sigma}}\Bigl(p_{\ka}e(\ka\cdot\tilde{a}_n\bal)\Bigr)\\
&+\biggl(\frac{1}{10\cdot2^{d}}\frac{M}{N}\sum_{\ka\in\mathbf{\Sigma}}(p_{\ka}e(\ka\cdot\tilde{a}_n\bal)\biggr)^2d\lambda(\bal)\Biggr)\\
\le&\prod_{n=1}^{\frac{N}{Q}}\Biggl(\biggl(1+\frac{1}{100\cdot2^{2d}}\biggl(\frac{M}{N}\biggr)^2\sum_{\ka\in\mathbf{\Sigma}}p_{\ka}^2\Biggr)\\
\le&\exp\left(\sum_{n=1}^{\frac{N}{Q}}\frac{1}{100}\sum_{\ka\in\mathbf{\Sigma}}\left|\restr{\hat{h}(\ka)}{\mathbf{z}_\taa=0}\right|^2\right)\\
\le&\exp\left(\frac{N}{Q}\frac{1}{100\cdot2^{d}}\left\|f\left(\V\frac{(\cdot)}{2}\right)\right\|_{L^2}^2+O\left(N^{-1000d}\right)\right)\\
=&\exp\left(\frac{N}{Q}\frac{1}{100\cdot2^{d}}\frac{M}{N}\left\|f\right\|_{L^2(\mathbb{T}^d)}^2+O\left(N^{-1000d}\right)\right)\le\exp\left(\frac{Q}{50\cdot2^{d}}\right).
\end{split}
\end{align}
\textbf{Step 3: Moment estimates (general measures):}
After defining $M$ and $Q$ as in \eqref{PQR}, and $p_{\ka}$ as in \eqref{pwiezuvor}, and changing the index of the dilated lacunary set of type $r$ to $\{\bal\ta_n\}_{K<n\le 2K}$, we look at the $d-$dimensional Fourier series of $$g(\bal):=\prod_{n=K+1}^{K+\frac{N}{Q}}\Biggl(1-\frac{1}{10\cdot2^{d}}\frac{M}{N}\sum_{\ka\in\mathbf{\Sigma}}\Bigl(p_{\ka}e(\ka\cdot\tilde{a}_n\bal)\Bigr)+\biggl(\frac{1}{10\cdot2^{d}}\frac{M}{N}\sum_{\ka\in\mathbf{\Sigma}}\bigl(p_{\ka}e(\ka\cdot\tilde{a}_n\bal)\bigr)\biggr)^2\Biggr).$$
$$g(\bal)=\sum_{\m\in\Z^d}q_{\m}e(\m\cdot\bal),$$
where
$$ q_{\m}=\int_{[0,1]^d}g(\bal)e(-\m\cdot \bal)d\lambda(\bal).$$
We apply the estimate \eqref{new1} and rewrite the first moment:
\begin{equation}
        \left\|e^{-\frac{1}{10}\omega^*_{N,\mathbf{Q}}}\right\|_{L^1(\mu)}\le\left\|g\right\|_{L^1(\mu)}=\left|\int_{[0,1]^d}g(\bal)d\mu(\bal)\right|
        \ll\sum_{\m\in\Z^d}\left|q_{\m}\right|\left(1+\|\m\|_{\ell^\infty}\right)^{-\Upsilon}, 
\end{equation}  
and split the sum into two parts, according to whether $\|\m\|_{\ell^\infty}\le a_N$ or $\|\m\|_{\ell^\infty}>a_N$.

The first sum can be treated identically to the calculations in \eqref{lebesguemeasure}--\eqref{lebesguemeasure2}, as the additional arising factors have frequencies too low compared to the existing ones in order to affect the calculation. This is again supported by \(L^2\)-orthogonality and yields the following bound:
\begin{align}\label{generalmeasure}
\begin{split}
   &\sum_{\|\m\|_{\ell^\infty}\le a_N} \left| q_{\m} \right| \left(1 +\| \m \|_{\ell^\infty} \right)^{-\Upsilon} \\
\le&\sum_{\|\m\|_{\ell^\infty}\le a_N}\int_{[0,1]^d}\prod_{n=K+1}^{K+\frac{N}{Q}}\Biggl(1-\frac{1}{10\cdot2^{d}}\frac{M}{N}\sum_{\ka\in\mathbf{\Sigma}}\Bigl(p_{\ka}e(\ka\cdot\tilde{a}_n\bal)\Bigr)\\
&+\biggl(\frac{1}{10\cdot2^{d}}\frac{M}{N}\sum_{\ka\in\mathbf{\Sigma}}\bigl(p_{\ka}e( \ka\cdot\tilde{a}_n\bal)\bigr)\biggr)^2\Biggr)e(-\mathbf{m}\cdot\bal)d\lambda(\bal)
\left(1 +\| \m \|_{\ell^\infty} \right)^{-\Upsilon}\\
\le&\exp\left(\sum_{n=1}^{\frac{N}{Q}}\frac{1}{100\cdot2^{d}}\frac{M}{N}\|f\|_{L^2(\mathbb{T}^d)}^2+O\left(N^{-1000d}\right)\right)\le\exp\left(\frac{Q}{50\cdot2^{d}}\right).
\end{split}
\end{align}
In the other sum, we make use of the rapid decay of the Fourier transform of the measure. We estimate as before:
\begin{align}\label{lebesguemeasure4}
\begin{split}
&\sum_{\|\m\|_{\ell^\infty}>a_N} \left| q_{\m} \right| \left(1 +\| \m \|_{\ell^\infty} \right)^{-\Upsilon} \\
   \le&\sum_{\|\m\|_{\ell^\infty}>a_N}\int_{[0,1]^d}\prod_{n=K+1}^{K+\frac{N}{Q}}\Biggl(1-\frac{1}{10\cdot2^{2d}}\frac{M}{N}\sum_{\ka\in\mathbf{\Sigma}}\Bigl(p_{\ka}e(\ka\cdot\tilde{a}_n\bal)\Bigr)\\
&+\biggl(\frac{1}{10\cdot2^{2d}}\frac{M}{N}\sum_{\ka\in\mathbf{\Sigma}}\bigl(p_{\ka}e(\ka\cdot\tilde{a}_n\bal)\bigr)\biggr)^2\Biggr)e( \m\cdot \bal)d\lambda(\bal)a_N^{-\Upsilon}.
\end{split}
\end{align}
Multiplying out the above product we get at most $3^{N/Q}$ terms, all of which have factors either $1$, $$\frac{1}{10R}\frac{M}{N}\sum_{\ka\in\mathbf{\Sigma}}p_\ka e(\ka\cdot\tilde{a}_n\bal),$$ 
for $n\in\{K+1,...,2N/Q\}$. Upon furhter multiplication, considering the sums in $\ka$, we get not more than $3^{N/Q}(2^{2d+1}N^2/M)^{2N/Q}=:(C(d)N^2/M)^{2N/Q}$ terms. Consequently, there are at most that many terms which are not orthogonal to some $m>a_N$.
This implies
\begin{align*}
   \begin{split}
 &\sum_{\|\m\|_{\ell^\infty}>a_N} \left| q_{\m} \right| \left(1 + \| \m \|_{\ell^\infty} \right)^{-\Upsilon} \le \left( \frac{C(d)N^2}{M} \right)^{\frac{2N}{Q}} a_N^{-\Upsilon}  \\
 =& \exp\left(\left((\log C(d))+2(\log N)-(\log M) \right)\frac{2N}{Q}-\Upsilon N   \right) <1\le\exp\left(\frac{Q}{50\cdot2^{d}}\right), 
\end{split} 
\end{align*}
for \(N\) large enough.\\


\textbf{Step 4: Almost sure existence:}
We need to cover $[0,1]^d$ by hyperrectangles $\mathbf{Q}\in\mathfrak{Q}$ which contain the support of $\omega_{N,\mathbf{Q}}$ for providing uniformity in $N$, i.e.\ there exists an $N_0\in\N$ such that for all $N\ge N_0$ every convex body $\mathcal{C}$ satisfies the statement of Theorem \ref{main}. Details will be explained in the next two steps. 

The calculation in step 2 and 3 have been carried out in order to allow an application of Markov's inequality and utilize the first Borel-Cantelli-Lemma. We have
\begin{align}\label{measure}
\begin{split}
    &\mu\left(\left\{\bal\in [0,1]^d\,:\,\frac{Q}{2^d}-\omega_{N,\mathbf{Q}}(\bal)\geq \frac{Q}{2^{d+1}}\right\}\right)\\\
    \le&\mu\left(\left\{\bal\in [0,1]^d\,:\,-\omega^{*}_{N,\mathbf{Q}}(\bal)\geq \frac{Q}{2^{d+2}}\right\}\right)\\\ 
    =&\mu\left(\left\{\bal\in [0,1]^d\,:\,\exp\left(-\frac{1}{10}\omega^{*}_{N,\mathbf{Q}}(\bal)\right)\geq \exp\left(\frac{Q}{40\cdot2^{d}}\right)\right\}\right)\\\
    \le& \frac{\|e^{-\frac{1}{10}\omega^*_{N,\mathbf{Q}}}\|_{L^1(\mu)}}{e^{\frac{Q}{40\cdot2^{d}}}}\le \frac{e^{\frac{Q}{50\cdot 2^{d}}}}{e^{\frac{Q}{40\cdot2^{d}}}}\le e^{-\frac{Q}{200\cdot 2^d}}=N^{-5d^2\max\{\Upsilon^{-1},1\}}\le N^{-5d^2}.
  \end{split}
\end{align}
Following this, we will need to consider at most $N^{3d^2}$ different $\mathbf{Q}\in\mathfrak{Q}$ regarding $\omega_{N,\mathbf{Q}}$, namely those $\mathbf{Q}\in\bar{\mathfrak{Q}}$ for $\bar{\mathfrak{Q}}$ defined below.

Since
\begin{equation}\label{bcant}
    \sum_{N\in\N}\sum_{\mathbf{Q}\in\overline{\mathfrak{Q}}}\mu(\{\bal\in [0,1]^d\,:\,\frac{Q}{2^d}-\omega_{N,\mathbf{Q}}(\bal)\geq \frac{Q}{2^{d+1}}\})\le\sum_{N\in\N}N^{3d^2}N^{-5d^2}<\infty,
\end{equation}
the first Borel Cantelli Lemma \eqref{bClambda} yields that for all sufficiently large $N$ and for almost all $\bal\in[0,1]^d$, we get
$$\omega_{N,\mathbf{Q}}(\bal)>0,$$
uniformly for all $\mathbf{Q}\in\bar{\mathfrak{Q}}$.
In other words, for every subset $\mathbf{Q}\in\bar{\mathfrak{Q}}$ there exists an $N_0$ such that for for all $N\ge N_0$ the hyperrectangle $\mathbf{Q}$ contains at least one point of $\{\tilde{a}_n\}_{n\le N}$.\\

\textbf{Step 5: Discretization:}
We will now construct a subset \(\overline{\mathfrak{Q}} \subseteq \mathfrak{Q}\) of cardinality at most \(N^{3d^2}\), ensuring that every convex body \(\mathcal{C}\) in \([0,1]^d\) with volume \(2^{d(\tilde{d}+2)}d^dM/N\) contains at least one element of \(\overline{\mathfrak{Q}}\).

We define
$$\overline{\mathfrak{Q}}=\left\{\mathbf{Q}\in\mathfrak{Q}|\mathbf{Q}=\mathbf{R}\mathbf{S}\cdot[-1,1]^d+\mathbf{t};\mathbf{R}\in\overline{\mathfrak{R}},\mathbf{S}\in\overline{\mathfrak{S}},\mathbf{t}\in\overline{\boldsymbol{\tau}}\right\},$$ 
where, for $1\le\delta<2$, that is allowed to depend on the matrix $\mathbf{S}$ below, we define 
\begin{align}
    \begin{split}
      \overline{\mathfrak{S}}&:=\left\{\mathbf{S}=\begin{pmatrix}
\frac{k_1}{N} & 0 & \cdots & 0 \\
0 & \frac{k_2}{N}  & \cdots & 0 \\
\vdots & \vdots  & \ddots & \vdots \\
0 & 0 &  \cdots & \frac{k_d\delta}{N}
\end{pmatrix}\middle|k_i\in[1,N]\cap\Z,\left\lceil\prod_{i=1}^{d}\frac{k_i}{N}\right\rceil=\frac{M}{2^dN},\delta \prod_{i=1}^{d}\frac{k_i}{N}=\frac{M}{2^dN} \right\},\\ 
 \overline{\boldsymbol{\tau}}&:=\left\{\mathbf{t}\in\R^d\middle|\mathbf{t}=\left(\frac{k_1}{N},\dots,\frac{k_d}{N}\right)^T,k_i\in[1,N]\cap\Z\right\}, \\
 \overline{\mathfrak{R}}&:=\left\{\mathbf{R}\in\mathfrak{R}\middle|\mathbf{R}=\prod_{l=1}^{\tilde{d}}\mathbf{R}_l\left(\frac{2\pi k_l}{N}\right),k_l\in[0,N-1]\cap\Z\right\}.\\
    \end{split}
\end{align}
Then $\overline{\mathfrak{Q}}$ is for sufficiently large $N$ of cardinality at most $N^d\cdot N^d\cdot N^{d(d-1)/2}\le N^{3d^2}$.

To achieve this, we
will prove the following three statements:
\begin{itemize}
    \item[Step 5.1:] Every hyperrectangle with origin-center of volume $2^dM/N$ that is not rotated nor translated (hence only a scaled hypercube) contains a hyperrectangle of volume $M/N$ constructed by scaling $[-1,1]^d$ by an element $\mathbf{S}\in\overline{\mathfrak{S}}$;   
    \item[Step 5.2:] Every rotated hyperrectangle with origin-center of volume $2^{d(\tilde{d}+1)}M/N$ contains a similar (proportional) hyperrectangle of volume $2^{d}M/N$, which is rotated by a rotation matrix $\mathbf{R}\in\overline{\mathfrak{R}}$ where we recall $\tilde{d}=d(d-1)/2$;
     \item[Step 5.3:] Every translated hyperrectangle of volume $2^{d(\tilde{d}+2)}M/N$ that is not rotated contains a homothetic hyperrectangle of volume $2^{d(\tilde{d}+1)}M/N$ translated by an element $\mathbf{t}\in\overline{\boldsymbol{\tau}}$.
\end{itemize}
Combining Theorem \ref{hard} with the above construction ensures that every convex body $\mathcal{C}$ of volume $2^{d(\tilde{d}+2)}d^dM/N$ contains an element from $\overline{\mathfrak{Q}}$.  This construction will be carried out in more detail below.

For every convex body \(\mathcal{C}\) of volume \(2^{d(\tilde{d}+2)}d^dM/N\), let us consider the homothetic hyperrectangles \(\mathfrak{T}_j\), \(j=1,\dots,\tilde{d}+2\), with the same center, where these hyperrectangles are constructed as follows: \(\mathfrak{T}_1\) is the hyperrectangle inscribed in \(\mathcal{C}\), as given by Theorem \ref{hard}, and they satisfy the volume relations
\begin{equation}\label{volumeassu}
   \frac{1}{d^d} \vol(\mathcal{C}) = \vol(\mathfrak{T}_1) = 2^{d(j-1)} \vol(\mathfrak{T}_j)= 2^{d(\tilde{d}+2)} \vol(\mathbf{Q})=2^{d(\tilde{d}+2)}\frac{M}{N}.
\end{equation}
Furthermore, we define 
\begin{equation}\label{defofD}
    D:=\vol(\mathfrak{T}_1).
\end{equation}
Observing definitions \eqref{matricesSRT}, we can write $\mathfrak{T}_j=\tilde{\mathbf{R}}\tilde{\mathbf{S}}_j\cdot[-1,1]^d+\tilde{\mathbf{t}}$, where $\tilde{\mathbf{S}}_j\in2^{\tilde{d}+3-j}\cdot\mathfrak{S}$ is a scaling matrix with $2^d$det$(\tilde{\mathbf{S}}_j)=\vol(\mathfrak{T_j})$, $\tilde{\mathbf{R}}\in\mathfrak{R}$ a rotation matrix, and $\tilde{\mathbf{t}}\in\boldsymbol{\tau}$ is a translation vector. 
In particular we have for $$\tilde{\mathbf{S}}_j=\begin{pmatrix}
\tilde{s}_{j,1} & 0  & \cdots & 0 \\
0 & \tilde{s}_{j,2}  & \cdots & 0 \\
\vdots & \vdots & \ddots & \vdots \\
0 & 0  & \cdots & \tilde{s}_{j,d}
\end{pmatrix},$$
that $\tilde{s}_{j,i}=2^{k-j}\tilde{s}_{k,i}$, for all $i=1,\dots,d$.

As elaborated upon earlier, we will divide this step into three smaller sub-steps, addressing scaling, rotation, and translation individually.

\textit{Step 5.1: Scaling:}
First we show that there exists an $\mathbf{S}\in\overline{\mathfrak{S}}$ such that:
\begin{equation}\label{step41claim}
    \mathbf{S}\cdot[-1,1]^d\subset\tilde{\mathbf{S}}_{\tilde{d}+2}\cdot[-1,1]^d.
\end{equation}
In particular we will prove the equivalent statement: There exists an $\mathbf{S}\in\overline{\mathfrak{S}}$ such that each diagonal matrix entry is smaller equal than the respective entry in $\tilde{\mathbf{S}}_{\tilde{d}+2}$, i.e.
we desire to find $k_i\in[1,N]\cap\Z,\ i=1,\dots,d-1$, such that
\begin{equation}\label{toshow1}
   \frac{k_i}{N}\le\tilde{s}_{\tilde{d}+2,i},\ \text{ and }\ \frac{k_d\delta}{N}\le\tilde{s}_{\tilde{d}+2,d}.
\end{equation}
The largest inscribed hyperrectangle in a convex body $\mathcal{C}\subset[0,1]^d$ have sidelength of at most $\sqrt{d}$, i.e.\ $s_{1,i}\le\sqrt{d}/2$. By construction \eqref{volumeassu} and the homothety assumption of $\mathfrak{T}_j$, we have for $j=1,\dots,\tilde{d}+2,\ \text{and}\ i=1,\dots,d$
\begin{equation}\label{upperbound}
    \tilde{s}_{j,i}\le \frac{\sqrt{d}}{2^{j}}.
\end{equation}
The above implies the following lower bound for every $s_{j,i}$:
\begin{equation}\label{sistklein}
    \frac{D}{2^j}=\frac{2^{jd}}{2^j}s_{j,i}\prod_{l=1,l\ne i}^d s_{j,l}\stackrel{\eqref{upperbound}}{\le}\frac{2^{jd}}{2^j}s_{j,i}\frac{\sqrt{d}}{2^{j(d-1)}}=s_{j,i}
\end{equation}
It is clear that we find for sufficiently large $N$ an $k_i\in[1,N]\cap\Z$, such that
\begin{equation}\label{clear}
    0<\tilde{s}_{\tilde{d}+2,i}-\frac{1}{N}\le\frac{k_i}{N}\le\tilde{s}_{\tilde{d}+2,i},
\end{equation}
for $i=1,\dots,d-1$, hence the first assertion in \eqref{toshow1} is proven. Moreover, this yields for sufficiently large $N$
\begin{equation}\label{smalltildes2}
   0\le \frac{D}{2^{\tilde{d}+3}}-\frac{1}{N}\stackrel{\eqref{sistklein}}{\le}\frac{\tilde{s}_{\tilde{d}+2,i}}{2}-\frac{1}{N}.
\end{equation}
Adding $s_{\tilde{d}+2,i}/2$ to the left and right handside of inequality \eqref{smalltildes2} together with \eqref{clear} implies
\begin{equation}\label{oneuse}
  \frac{\tilde{s}_{\tilde{d}+2,i}}{2}< \tilde{s}_{\tilde{d}+2,i}-\frac{1}{N}<\frac{k_i}{N}., 
\end{equation}
Utilizing the above equation \eqref{oneuse} yields
\begin{equation}
    \frac{k_d\delta}{N}=\frac{\frac{D}{2^{(\tilde{d}+3)d}}}{\prod_{i=1}^{d-1}\frac{k_i}{N}}\stackrel{\eqref{oneuse}}{<}\frac{\frac{D}{2^{(\tilde{d}+3)d}}}{\prod_{i=1}^{d-1}\frac{\tilde{s}_{\tilde{d}+2,i}}{2}}=\frac{\tilde{s}_{\tilde{d}+2,d}}{2},
\end{equation}
which implies the second inequality in \eqref{toshow1}, hence the claim \eqref{step41claim}.

\begin{figure}[h]
    \centering
    \begin{tikzpicture}

    \draw[blue] (-4, -2) rectangle (4, 2);
    \node[blue] at (-2.85,1.65){$\tilde{\mathbf{S}}_{\tilde{d}+2}
    \cdot[-1,1]^d$};
    \node[blue][below] at (0,-2){$2\tilde{s}_{\tilde{d}+2,1}$};
    \node[blue][right] at (4,0){$2\tilde{s}_{\tilde{d}+2,2}$};

    \draw[red] (-3.5, -0.57142857142) rectangle (3.5, 0.57142857142);
    \node[red] at (-2.5, 0.22142857142){$\mathbf{S}\cdot[-1,1]^d$};
    \node[red][below] at (0, -0.57142857142){$2s_{\tilde{d}+2,1}$};
    \node[red][left] at (3.5, 0){$2s_{\tilde{d}+2,2}$};
    \end{tikzpicture}
    
    \caption{$\mathbf{S}\cdot[-1,1]^d\subset\tilde{\mathbf{S}}_{\tilde{d}+2}\cdot[-1,1]^d$ for $d=2$.}
    \label{fig:enter-label}
\end{figure}
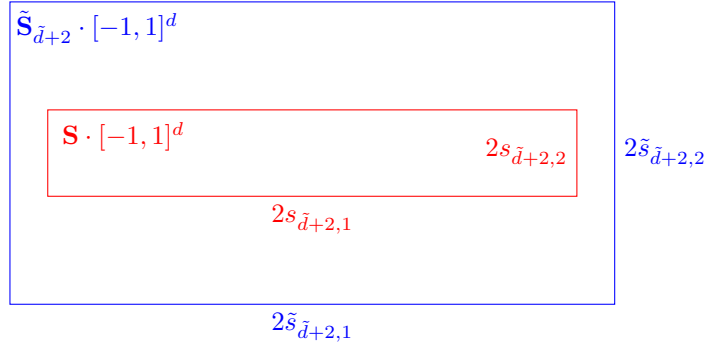

\textit{Step 5.2: Rotation:} Next we show that there exists an $\mathbf{R}\in\overline{\mathfrak{R}}$ such that 
\begin{equation}\label{rot}
\mathbf{R}\tilde{\mathbf{S}}_{\tilde{d}+2}\cdot[-1,1]^d\subset\tilde{\mathbf{R}}\tilde{\mathbf{S}}_2\cdot[-1,1]^d.
\end{equation}
As we will see at the end of this sub-section, it is sufficient to prove 
\begin{equation}\label{rotsimple}
\mathbf{R}_n(2\pi/N)\tilde{\mathbf{S}}_{j}\cdot[-1,1]^d\subset\tilde{\mathbf{S}}_{j+1}\cdot[-1,1]^d,
\end{equation}
for any $n=1,\dots,\tilde{d}$ and any $j=2,\dots,\tilde{d}+2$.
Since the procedure is the same for every $j$ and $n$, 
we will prove inclusion \eqref{rotsimple} only for $n=1$ and without loss of generality let us omit the index $j$ until it is needed again.
Thus, we write $\tilde{\mathbf{S}}$ instead of $\tilde{\mathbf{S}}_j$ and $\tilde{s}_i:=\tilde{s}_{j,i}$, $i=1,2$. 

We rotate the vectors $\mathbf{v}_1=(\tilde{s}_1,-\tilde{s}_2)^T$ and $\mathbf{v}_2=(\tilde{s}_1,\tilde{s}_2)^T$ by the angle $2\pi/N$, considering only the first entry of $\mathbf{v_1}$ and the second of $\mathbf{v_2}$, and subtract $\tilde{s}_1$, $\tilde{s}_2$ respectively, see Figure \ref{rotation}.
The larger of the above two is equal to the Hausdorff distance between the sets $\mathbf{R}_1(2\pi/N)\tilde{\mathbf{S}}\cdot[-1,1]^d$ and $\tilde{\mathbf{S}}\cdot[-1,1]^d$. These two expressions in the first and second dimensional direction correspond to the following respectively:
\begin{subequations} \label{greatdist}
 \begin{align}
\text{dist}_1\left(\frac{2\pi}{N}\right)&:=\tilde{s}_1\cos\left(\frac{2\pi}{N}\right)+\tilde{s}_2\sin\left(\frac{2\pi}{N}\right)-\tilde{s}_1,\\
\text{dist}_2\left(\frac{2\pi}{N}\right)&:=\tilde{s}_1\sin\left(\frac{2\pi}{N}\right)+\tilde{s}_2\cos\left(\frac{2\pi}{N}\right)-\tilde{s}_2.
 \end{align}
 \end{subequations}
For 
\begin{equation}\label{rotsimple3}
\mathbf{R}_1\left(\frac{2\pi}{N}\right)\tilde{\mathbf{S}}_{j+1}\cdot[-1,1]^d\subset\tilde{\mathbf{S}}_{j}\cdot[-1,1]^d
\end{equation}
to be true we need show 
\begin{equation} \label{greatdist2}
\text{dist}_1\left(\frac{2\pi}{N}\right)\le\tilde{s}_1,\quad\text{dist}_2\left(\frac{2\pi}{N}\right)\le\tilde{s}_2.
 \end{equation}
The second assertion follows for $N$ sufficiently large directly from 
\begin{equation}
    \tilde{s}_1\sin\left(\frac{2\pi}{N}\right)<\frac{1}{2}\sin\left(\frac{2\pi}{N}\right)<\frac{D}{2^{\tilde{d}+2}}<\tilde{s}_{2}<\tilde{s}_{2}\left(2-\cos\left(\frac{2\pi}{N}\right)\right),
\end{equation}
where we employed the properties \eqref{upperbound} and \eqref{sistklein}. The first assertion can be shown identically.
Moreover, since for $N$ sufficiently large 
\begin{equation}
    \text{dist}_i(\theta)\le\text{dist}_i\left(\frac{2\pi}{N}\right),\quad i=1,2
\end{equation}
is true for every angle $0\le\theta\le2\pi/N$, we have 
\begin{equation}\label{rotsimple4}
\mathbf{R}_1(\theta)\tilde{\mathbf{S}}_{j+1}\cdot[-1,1]^d\subset\tilde{\mathbf{S}}_{j}\cdot[-1,1]^d.
\end{equation}
Now we observe, that for any $\tilde{\gamma}\in[0,2\pi]$, there exists a $k\in[0,N-1]\cap\Z$, such that
\begin{equation}\label{ungleichunga2}
    \tilde{\gamma}\le\frac{2\pi k}{N}\le\tilde{\gamma}+\frac{2\pi}{N},
\end{equation}
which, together with inclusion \eqref{rotsimple4}, implies
for $\tilde{\mathbf{R}}=\prod_{l=0}^{\tilde{d}-1}\tilde{\mathbf{R}}_{\tilde{d}-l}(\tilde{\gamma}_{\tilde{d}-l})$, there exists $k_n$ with $n=1,\dots,\tilde{d}$ such that

\begin{equation}\label{rotsimple7}
\mathbf{R}_{n}\left(\frac{k_n2\pi}{N}\right)\tilde{\mathbf{S}}_{j}\cdot[-1,1]^d\subset\tilde{\mathbf{R}}_{n}(\tilde{\gamma_n})\tilde{\mathbf{S}}_{j+1}\cdot[-1,1]^d.
\end{equation}
Thus there exists an $k_l\in[0,N-1]\cap\Z$, for all $l=1,\dots,\tilde{d}$, such that for $\mathbf{R}=\prod_{l=0}^{\tilde{d}-1}\mathbf{R}_{\tilde{d}-l}(k_{\tilde{d}-l}2\pi/N)\in\mathfrak{R}$ we have
\begin{align}
\begin{split}
        &\mathbf{R}\tils_{\tilde{d}+2}[-1,1]^d=\ra_{\tilde{d}}\left(\frac{k_{\tilde{d}}2\pi}{N}\right)\cdots\ra_1\left(\frac{k_12\pi}{N}\right)\tils_{\tilde{d}+2}[-1,1]^d\\
        \subset&\ra_{\tilde{d}}\left(\frac{k_{\tilde{d}}2\pi}{N}\right)\cdots\ra_{2}\left(\frac{k_22\pi}{N}\right)\tilr_{1}\left(\tilga_1\right)\tils_{\tilde{d}+1}[-1,1]^d\\
        \subset&\cdots\subset\tilr_{\tilde{d}}\left(\tilga_{\tilde{d}}\right)\dots\tilr_{1}\left(\tilga_{1}\right)\tils_{2}[-1,1]^d=\tilr\tils_{2}[-1,1]^d.
\end{split}
\end{align}

\begin{figure}[h]
    \centering
    \begin{tikzpicture}

    \draw[blue] (-4, -2) rectangle (4, 2);
    \node[blue][right] at (-4,1.67){$\tilde{\mathbf{R}}\tilde{\mathbf{S}}_{j-1}
    \cdot[-1,1]^d$};
    \node[blue][below] at (0,-2){$2\tilde{s}_{j-1,1}$};
    \node[blue][right] at (4,0){$2\tilde{s}_{j-1,2}$};
    \draw[green] (-2, -1) rectangle (2, 1);
    \node[green][right] at (-2, 1.25){$\tilde{\mathbf{R}}_n\tilde{\mathbf{S}}_j\cdot[-1,1]^d$};
    \node[green][below] at (0, -1.1){$2\tilde{s}_{j,1}$};
    \node[green][right] at (2, 0){$2\tilde{s}_{j,2}$};
    \draw[red] (-1.6730326075, -1.48356) -- (-2.1906706977, 0.4482877361);
    \draw[red] (-2.1906706977, 0.4482877361) -- (1.6730326075, 1.48356);
    \draw[red] (-1.6730326075, -1.48356) -- (2.1906706977, -0.4482877361);
    \draw[red] (2.1906706977, -0.4482877361) -- (1.6730326075, 1.48356);
    \node[red][right] at (-2, 0.15){$\mathbf{R}\tilde{\mathbf{S}}_j\cdot[-1,1]^d$};
    \draw[purple] (2.1906706977, -0.4482877361) -- (2, -0.4482877361);
    \node[purple][below] at (2.3906706977, -0.4482877361) {$\text{dist}_{1}(\theta)$};
    \draw[purple] (1.6730326075, 1.48356) -- (1.6730326075, 1);
    \node[purple][right] at (1.7730326075, 1.28356) {$\text{dist}_{2}(\theta)$};
    \draw[brown] (-1.5373650982, -1.6237329074) -- (-2.2214053849, 0.25565);
    \draw[brown] (-2.2214053849, 0.25565) -- (1.5373650982, 1.6237329074);
    \draw[brown] (-1.5373650982, -1.6237329074) -- (2.2214053849, -0.25565);
    \draw[brown] (2.2214053849, -0.25565) -- (1.5373650982, 1.6237329074);
    \node[brown][above] at (-0.25, -0.8){$\mathbf{R}_1\left(\frac{2\pi}{N}\right)\tilde{\mathbf{R}}\tilde{\mathbf{S}}_j\cdot[-1,1]^d$};
    
    \end{tikzpicture}
    
    \caption{$\mathbf{R}\tilde{\mathbf{S}}_3\cdot[-1,1]^d\subset\tilde{\mathbf{R}}\tilde{\mathbf{S}}_2\cdot[-1,1]^d$ for $d=2$.}
    \label{rotation}
\end{figure}
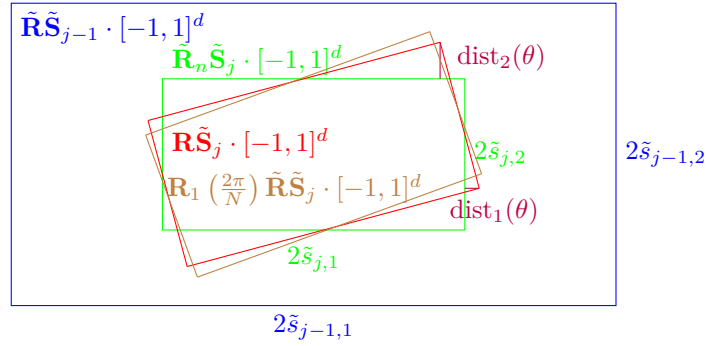

\textit{Step 5.3: Translation:}
Last, we show that there exists an $\mathbf{t}\in\overline{\boldsymbol{\tau}}$ such that 
\begin{equation}\label{transl}   \tilde{\mathbf{S}}_2\cdot[-1,1]^d+\taa\subset\tilde{\mathbf{S}}_1\cdot[-1,1]^d+\tilde{\taa}:
\end{equation}
By construction, it is sufficient to show that we can find an $\taa\in\overline{\boldsymbol{\tau}}$ such that $\tilde{\mathbf{S}}_2\cdot[-1,1]^d+\taa$ is $\tilde{\mathbf{S}}_2\cdot[-1,1]^d+\tilde{\taa}$ shifted in any dimensional direction by not more than $\tilde{s}_{2,i}$, $i=1,\dots,d$ respectively. 

We observe, that there exist $k_i\in[1,N]\cap\Z$, such that for $\tilde{\mathbf{t}}=\left(\tilde{t}_1,\dots,\tilde{t}_d\right)^T$ we have
\begin{equation}\label{transineq}
    \tilde{t}_i-\frac{1}{N}\le\frac{k_i}{N}\le\tilde{t}_i,
\end{equation}
for every $i=1,\dots,d$, and for $N$ sufficiently large, we can estimate
\begin{equation}\label{transineq2}
    \frac{1}{N}<\frac{D}{4}<\tilde{s}_{2}.
\end{equation}
The inequality \eqref{transineq} implies that there exists an $\taa\in\overline{\tau}$ such that $\tilde{\mathbf{S}}_2\cdot[-1,1]^d+\taa$ is just $\tilde{\mathbf{S}}_2\cdot[-1,1]^d+\tilde{\taa}$ translated by a vector with entries at most $1/N$, which together with \eqref{transineq2} implies the desired claim \eqref{transl}.

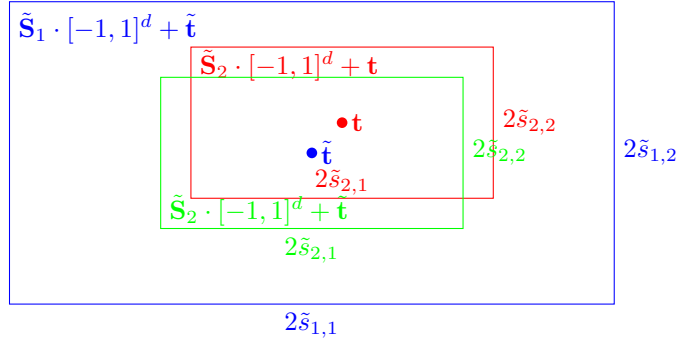
\begin{figure}[h]
    \centering
    \begin{tikzpicture}

    \draw[blue] (-4, -2) rectangle (4, 2);
    \node[blue] at (-2.7,1.65){$\tilde{\mathbf{S}}_1
    \cdot[-1,1]^d+\tilde{\mathbf{t}}$};
    \node[blue][below] at (0,-2){$2\tilde{s}_{1,1}$};
    \node[blue][right] at (4,0){$2\tilde{s}_{1,2}$};
    \fill[blue] (0,0) circle[radius=2pt];
    \node[blue][right] at (0,0){$\tilde{\taa}$};
    \draw[red] (-1.6, -0.6) rectangle (2.4, 1.4);
    \node[red] at (-0.3, 1.16){$\tilde{\mathbf{S}}_2\cdot[-1,1]^d+\mathbf{t}$};
    \node[red][below] at (0.4, -0.1){$2\tilde{s}_{2,1}$};
    \node[red][right] at (2.4, 0.4){$2\tilde{s}_{2,2}$};
    \fill[red] (0.4,0.4) circle[radius=2pt];
    \node[red][right] at (0.4,0.4){$\taa$};
    \draw[green] (-2, -1) rectangle (2, 1);
    \node[green] at (-0.7, -0.76){$\tilde{\mathbf{S}}_2\cdot[-1,1]^d+\tilde{\mathbf{t}}$};
    \node[green][below] at (0, -1){$2\tilde{s}_{2,1}$};
    \node[green][right] at (2, 0){$2\tilde{s}_{2,2}$};
    \end{tikzpicture}
    
    \caption{$\tilde{\mathbf{S}}_{2}\cdot[-1,1]^d+\mathbf{t}\subset\tilde{\mathbf{S}}_1\cdot[-1,1]^d+\tilde{\taa}$ for $d=2$.}
    \label{fig:enter-label2}
\end{figure}

By combining all the steps, we get for every convex body $\mathcal{C}$, that there exists a $\mathbf{Q} \in \overline{\mathfrak{Q}}$ such that
\begin{align}\label{inclus}
    \begin{split}
        \mathbf{Q} &= \mathbf{R}\mathbf{S} \cdot [-1, 1]^d + \mathbf{t} 
\subset \mathbf{R}\tilde{\mathbf{S}}_{\tilde{d}+2} \cdot [-1, 1]^d + \taa\\
&\subset \tilde{\mathbf{R}}\tilde{\mathbf{S}}_2 \cdot [-1, 1]^d + \taa\subset \tilde{\mathbf{R}}\tilde{\mathbf{S}}_1 \cdot [-1, 1]^d + \tilde{\taa}
 = \mathfrak{T}_1 \subset \mathcal{C},
    \end{split}
\end{align}
with $\vol(\mathcal{C}) \leq 2^{d(\tilde{d}+2)}d^d \vol(\mathbf{Q})$, and $|\overline{\mathfrak{Q}}| \le N^{3d^2}$.\\


\textbf{Step 6: Convex sets:}
From Step 5, we conclude that for \(\mu\)-almost every \(\boldsymbol{\alpha} \in [0,1]^d\), there exists an \(N_0 \in \mathbb{N}\) such that, for every \(N \geq N_0\), each \(\mathbf{Q} \in \overline{\mathfrak{Q}}\) contains at least one point of the truncated set \(\{\boldsymbol{\alpha} a_n\}_{n \leq N}\).

By the previous step, we know that for any convex body \(\mathcal{C}\) with volume \(2^{d(\tilde{d}+2)}d^d(\log N)^2/N\), there exists a hyperrectangle \(\mathbf{Q} \in \overline{\mathfrak{Q}}\) that is entirely contained within \(\mathcal{C}\).

Since
\[
\vol(\mathcal{C}) = 2^{d(\tilde{d}+2)}d^d \frac{(\log N)^{2}}{N} \simeq \frac{(\log N)^{2}}{N},
\]
which proves the theorem.


\end{proof}

\section{\textbf{Appendix:}\label{Appendix} Proof of Hadwiger's Theorem \ref{hard}}

For the convenience of the reader, in this section we provide a translated version of Hadwiger's proof \cite{378850199_0010}. 


\begin{figure}[h]
    \centering
\resizebox{1\textwidth}{!}{%
\begin{circuitikz}
\tikzstyle{every node}=[font=\LARGE]
\draw [ line width=1pt ] (4,13.25) rectangle (17,8);
\draw [short] (17,10) -- (4,10);
\draw [ color={rgb,255:red,26; green,95; blue,180}, short] (17,10) -- (16.75,11.25);
\draw [ color={rgb,255:red,26; green,95; blue,180}, short] (16.75,11.25) -- (15.75,12.25);
\draw [ color={rgb,255:red,26; green,95; blue,180}, short] (15.75,12.25) -- (15,12.5);
\draw [ color={rgb,255:red,26; green,95; blue,180}, short] (7,13.25) -- (6.25,13);
\draw [ color={rgb,255:red,26; green,95; blue,180}, short] (6.25,13) to[out=205,in=80] (4,10);
\draw [ color={rgb,255:red,26; green,95; blue,180}, short] (4,10) -- (6,8);
\draw [ color={rgb,255:red,26; green,95; blue,180}, short] (6,8) -- (11.75,8.25);
\draw [ color={rgb,255:red,26; green,95; blue,180}, short] (11.75,8.25) to[out=5,in=230] (17,10);
\node [font=\LARGE] at (3.6,10.2) {H};
\draw [short] (17,10) -- (11.25,8.25);
\draw [short] (11.25,8.25) -- (4,10);
\draw [short] (11.25,8) -- (11.25,13.25);
\draw [short] (4,10) -- (11.25,13);
\draw [short] (11.25,13) -- (17,10);
\draw [, line width=1pt ] (7.5,11.41) rectangle (14.25,9.19);
\draw [ color={rgb,255:red,26; green,95; blue,180}, short] (7,13.25) -- (11.25,13);
\draw [ color={rgb,255:red,26; green,95; blue,180}, short] (11.25,13) -- (15,12.5);
\draw [short] (7.5,13.25) -- (7.5,8);
\draw [short] (14.25,13.25) -- (14.25,8);
\draw [short] (10.5,13.25) -- (10.5,8);
\node [font=\LARGE] at (10.5,13.75) {0};
\node [font=\LARGE] at (11.25,13.75) {s};
\node [font=\large] at (7.5,13.75) {1/2(s+h)};
\node [font=\large] at (14.25,13.75) {1/2(s-h)};
\node [font=\large] at (4,13.75) {h};
\node [font=\large] at (17,13.75) {-h};
\node [font=\LARGE] at (17.4,10.35) {  H'};
\draw [->, >=Stealth] (4,13.25) -- (4,13.5);
\draw [->, >=Stealth] (7.5,13.25) -- (7.5,13.5);
\draw [->, >=Stealth] (10.5,13.25) -- (10.5,13.5);
\draw [->, >=Stealth] (11.25,13.25) -- (11.25,13.5);
\draw [->, >=Stealth] (14.25,13.25) -- (14.25,13.5);
\draw [->, >=Stealth] (17,13.25) -- (17,13.5);
\draw [->, >=Stealth] (17,10) -- (18,10);
\node [font=\Large] at (18.25,10) {z};
\end{circuitikz}
}%
    
    \caption{}
    \label{Hadwiger}
\end{figure}

\begin{proof}[Proof of Theorem \eqref{hard}]
    
The claims \eqref{a} and \eqref{b} are correct and trivial in the case \( d = 1 \). Let \( k > 1 \), and we make the inductive assumption that our statements are already established for dimension \( d - 1 \).

Consider now the \( d \)-dimensional case. We start from a diameter chord \( HH' \) of the convex body \( \mathcal{C} \) of length \( 2h \). Such a chord is the connection of two points \( H \) and \( H' \), which realize the greatest distance provided by the point pairs of the convex body. Compare Figure \ref{Hadwiger}, which refers to the case \( d = 2 \). We choose the middle of the chord \( HH' \) as the origin of a coordinate system and let the chord fall into what we call the $z$-axis. Let \( E(t) \) denote the \((d-1)\)-dimensional plane \( z = t \), orthogonal to the $z$-axis. By construction, \( E(h) \) and \( E(-h) \) are parallel support planes of \( \mathcal{C} \).
Indeed, the convex body must be entirely in the closed parallel strip between \( E(h) \) and \( E(-h) \), as any point outside these planes would exceed the value \( 2h \), contradicting the definition of the diameter.

Through orthogonal projection from \( \mathcal{C} \) to \( E(0) \), we obtain a \((d-1)\)-dimensional convex body \( \mathcal{C}_0 \) of volume \( \vol(\mathcal{C}_0) \). According to the inductive assumption, there is a hyperrectangle \( P_0 \) in \( E(0) \), such that \( \mathcal{C}_0 \subseteq P_0 \) and the inequality \begin{equation}\label{aa}
\vol(P_0) \leq (k - 1)! \vol(\mathcal{C}_0)
\end{equation}
holds. On the other hand, the relation
\begin{equation}\label{ab}
2h \vol(\mathcal{C}_0) \leq k \vol(\mathcal{C})
\end{equation}
also holds.

 To see this, we first consider the body \( \overline{\mathcal{C}} \) symmetrical about \( E(0) \), which arises from \( \mathcal{C} \) by a Steiner symmetrization, on the plane \( E(0) \), see e.g. \cite{gardner1995geometric} for details. From the simplest properties of this transformation it follows that $H,H'\subset\overline{\mathcal{C}}$ and \( \mathcal{C}_0 \subseteq \overline{\mathcal{C}} \). We form the convex hull of \( H \), \( H' \), and \( \mathcal{C}_0 \), thus creating a symmetrical double cone body \( X \). Hence, \( X \subseteq \overline{\mathcal{C}} \), so \( \vol(X) \leq \vol(\overline{\mathcal{C}}) \). Since further \( \vol(\mathcal{C}) = \vol(\overline{\mathcal{C}}) \), equation \eqref{ab} follows.
Based on the construction and the given elementary geometric situation, it is now easy to see that \( \mathcal{C} \) is covered by a hyperrectangle \( P \), whose cross-sectional hyperrectangle in \( E(0) \) is identical to \( P_0 \), with height \( 2h \) and volume
\begin{equation}\label{ac}
    \vol(P) = 2h \vol(P_0) \quad [\mathcal{C} \subseteq P].
\end{equation}
The combination of the relations \eqref{aa}, \eqref{ab}, and \eqref{ac} proves assertion \eqref{a}.

We choose a new starting point with the volume formula
\[
\vol(\mathcal{C}) = \int_{-h}^{h} \vol(\mathcal{C} \cap E(t)) \, dt,
\]
where we integrate over the continuously varying \((d-1)\)-dimensional volumes of the cross-sections \( \mathcal{C} \cap E(t) \).
According to the mean value theorem, there exists an \( s \), with \( -h < s < h \) such that with the setting \( \mathcal{C} \cap E(s) = \mathcal{C}_{00} \) the relation
\begin{equation}\label{ba}
    2h \vol(\mathcal{C}_{00}) = \vol(A)
\end{equation}
holds. According to the inductive assumption, there is a hyperrectangle \( Q_0 \) in the space \( E(s) \) so that \( Q_0 \subseteq \mathcal{C}_{00} \) holds and the inequality
\begin{equation}\label{bb}
    \vol(Q_0) \geq \left( \frac{1}{d-1} \right)^{d-1} \vol(\mathcal{C}_{00})
\end{equation}
holds. If we now form the convex hull of \( H, H' \), and \( Q_0 \), a symmetrical double pyramid \( Y \) is created, and \( Y \subseteq A \). If \( 0 < \lambda < 1 \), then the sections
\[
S_{\lambda} = Y \cap E([(1-\lambda)s + \lambda h]),
\]
\[
S'_{\lambda} = Y \cap E([(1-\lambda)s - \lambda h])
\]
are two translationally equivalent hyperrectangles in the $z$-direction, which arise from \( Q_0 \) by similar scaling in the ratio \( 1 : (1 - \lambda) \). \( S_{\lambda} \) and \( S'_{\lambda} \) form the base and top surfaces of a hyperrectangle \( Q_{\lambda} \) of height \( 2\lambda h \), which, as is clear from the construction, is contained in \( A \), and whose volume is given by

\begin{equation}\label{bc}
    \vol(Q_{\lambda}) = 2 \lambda h (1 - \lambda)^{d-1} \vol(Q_0).
\end{equation}

If we set \(\lambda = 1/d\) and \(Q_{\lambda} = Q\), then from equations \eqref{ba}, \eqref{bb}, and \eqref{bc} the assertion \eqref{b} follows.

\end{proof} 

\textbf{Acknowledgments}
This research was funded in whole, or in part, by the Austrian Science Fund (FWF)  [Grant-DOI 10.55776/P35322]. 
The author expresses gratitude to Christoph Aistleitner, Athanasios Sourmelidis, Andrei Shubin and Maryna Manskova for many helpful discussions. He is also grateful to Faustin Adiceam for pointing out the connection with Danzer's problem, and to Niclas Technau, Matteo Bordignon, and Giuseppe Molteni for helpful feedback. The author is also indebted to the anonymous reviewer for valuable corrections and helpful suggestions that improved the manuscript.

\bibliographystyle{siam}
\bibliography{main}

\begin{thebibliography}{10}

\bibitem{adiceam}
{\sc F.~Adiceam}, {\em Around the {D}anzer problem and the construction of dense forests}, Enseign. Math., 68 (2022), pp.~25--60.

\bibitem{abt}
{\sc C.~Aistleitner, I.~Berkes, and R.~Tichy}, {\em Lacunary sequences in analysis, probability and number theory}, in Diophantine problems: determinism, randomness and applications, vol.~62 of Panor. Synth\`eses, Soc. Math. France, Paris, 2024, pp.~1--60.

\bibitem{ahr}
{\sc C.~Aistleitner, A.~Hinrichs, and D.~Rudolf}, {\em On the size of the largest empty box amidst a point set}, Discrete Appl. Math., 230 (2017), pp.~146--150.

\bibitem{al}
{\sc A.~Arman and A.~E. Litvak}, {\em Minimal dispersion on the cube and the torus}, J. Complexity, 85 (2024), pp.~Paper No. 101883, 7.

\bibitem{ball}
{\sc K.~Ball}, {\em An elementary introduction to modern convex geometry}, in Flavors of geometry, vol.~31 of Math. Sci. Res. Inst. Publ., Cambridge Univ. Press, Cambridge, 1997, pp.~1--58.

\bibitem{bugi}
{\sc Y.~Bugeaud}, {\em Around the littlewood conjecture in diophantine approximation}, in Numéro consacré au trimestre ``Méthodes arithmétiques et applications'', automne 2013, vol.~2014/1 of Publ. Math. Besançon Algèbre Théorie Nr., Presses Univ. Franche-Comté, Besançon, 2014, pp.~5--18.

\bibitem{chye}
{\sc S.~Chaubey and N.~Yesha}, {\em The distribution of spacings of real-valued lacunary sequences modulo one}, Mathematika, 68 (2022), pp.~416--428.

\bibitem{chow2023dispersion}
{\sc S.~Chow and N.~Technau}, {\em Dispersion and littlewood's conjecture}, Advances in Mathematics, 447 (2024), p.~109697.

\bibitem{demeter2020fourier}
{\sc C.~Demeter}, {\em Fourier Restriction, Decoupling and Applications}, Cambridge Studies in Advanced Mathematics, Cambridge University Press, 2020.

\bibitem{dts}
{\sc M.~Drmota and R.~F. Tichy}, {\em Sequences, discrepancies and applications}, vol.~1651 of Lecture Notes in Mathematics, Springer-Verlag, Berlin, 1997.

\bibitem{dj}
{\sc A.~Dumitrescu and M.~Jiang}, {\em On the largest empty axis-parallel box amidst {$n$} points}, Algorithmica, 66 (2013), pp.~225--248.

\bibitem{durrett2019probability}
{\sc R.~Durrett and R.~Durrett}, {\em Probability: Theory and Examples}, Cambridge Series in Statistical and Probabilistic Mathematics, Cambridge University Press, 2019.

\bibitem{fuku}
{\sc K.~Fukuyama}, {\em The law of the iterated logarithm for discrepancies of {$\{\theta^nx\}$}}, Acta Math. Hungar., 118 (2008), pp.~155--170.

\bibitem{gardner1995geometric}
{\sc R.~J. Gardner}, {\em Geometric Tomography}, Encyclopedia of Mathematics and its Applications, Cambridge University Press, 1995.

\bibitem{gröchenig2001foundations}
{\sc K.~Gr{\"o}chenig}, {\em Foundations of Time-Frequency Analysis}, Applied and Numerical Harmonic Analysis, Birkh{\"a}user Boston, 2001.

\bibitem{378850199_0010}
{\sc H.~Hadwiger}, {\em Volumschätzung für die einen Eikörper überdeckenden und unterdeckenden Parallelotope}, Birkhäuser, 1955.

\bibitem{hkkr}
{\sc A.~Hinrichs, D.~Krieg, R.~J. Kunsch, and D.~Rudolf}, {\em Expected dispersion of uniformly distributed points}, J. Complexity, 61 (2020), pp.~101483, 9.

\bibitem{john}
{\sc F.~John}, {\em Extremum problems with inequalities as subsidiary conditions}, in Studies and {E}ssays {P}resented to {R}. {C}ourant on his 60th {B}irthday, {J}anuary 8, 1948, Interscience Publishers, New York, 1948, pp.~187--204.

\bibitem{litvak}
{\sc A.~E. Litvak}, {\em A remark on the minimal dispersion}, Commun. Contemp. Math., 23 (2021), pp.~Paper No. 2050060, 18.

\bibitem{phil}
{\sc W.~Philipp}, {\em Limit theorems for lacunary series and uniform distribution {${\rm mod}\ 1$}}, Acta Arith., 26 (1974/75), pp.~241--251.

\bibitem{rotetichy}
{\sc G.~Rote and R.~F. Tichy}, {\em Quasi-{M}onte {C}arlo methods and the dispersion of point sequences}, vol.~23, 1996, pp.~9--23.
\newblock Monte Carlo and quasi-Monte Carlo methods.

\bibitem{rudza}
{\sc Z.~Rudnick and A.~Zaharescu}, {\em The distribution of spacings between fractional parts of lacunary sequences}, Forum Math., 14 (2002), pp.~691--712.

\bibitem{sw}
{\sc Y.~Solomon and B.~Weiss}, {\em Dense forests and {D}anzer sets}, Ann. Sci. \'{E}c. Norm. Sup\'{e}r. (4), 49 (2016), pp.~1053--1074.

\bibitem{stefanescu2024dispersiondilatedlacunarysequences}
{\sc E.~Stefanescu}, {\em The dispersion of dilated lacunary sequences, with applications in multiplicative diophantine approximation}, Advances in Mathematics, 461 (2025), p.~110062.

\bibitem{fourier}
{\sc E.~M. Stein and R.~Shakarchi}, {\em Fourier Analysis: An Introduction}, Princeton University Press, 2003.

\bibitem{tz}
{\sc N.~Technau and A.~Zafeiropoulos}, {\em The discrepancy of {$(n_kx)^\infty_{k=1}$} with respect to certain probability measures}, Q. J. Math., 71 (2020), pp.~573--597.

\bibitem{tem}
{\sc V.~N. Temlyakov}, {\em Universal discretization}, J. Complexity, 47 (2018), pp.~97--109.

\bibitem{tvv}
{\sc M.~Tr\"{o}dler, J.~Volec, and J.~Vybiral}, {\em A tight lower bound on the minimal dispersion}, European J. Combin., 120 (2024), pp.~Paper No. 103945, 5.

\bibitem{yesha}
{\sc N.~Yesha}, {\em Intermediate-scale statistics for real-valued lacunary sequences}, Math. Proc. Cambridge Philos. Soc., 175 (2023), pp.~303--318.

\end{thebibliography}







\end{document}